%
%
%
 \magnification \magstep1
\font\bold=cmbx10 at 12pt

\font\tts=cmtt8
\centerline{\bold A commutative noetherian local ring of embedding dimension 4 having a }
\bigskip
\centerline{\bold transcendental series of Betti numbers}
\bigskip
\centerline{ Jan-Erik Roos}
\centerline{Department of Mathematics}
\centerline{Stockholm University}
\centerline{ e-mail: {\tt  jeroos@math.su.se}}
\bigskip
\centerline{February 14, 2017}

\bigskip

\def\mysec#1{\bigskip\centerline{\bf #1}\nobreak\par}

\def\cite#1{~[{\bf #1}]}
{ \openup 1\jot
\mysec{ Abstract}
We construct a ring with the properties of the title of the paper.
We also construct some other local rings of embedding dimension 4
with exotic properties. Among the methods used are the {\tt Macaulay2}-package {\tt DGAlgebras}
 by Frank Moore, combined with and inspired by results by Anick, Avramov, Backelin, Katth\"an, Lemaire,
Levin, L\"ofwall and others.
\medskip
{\it Mathematics Subject Classification (2010):} Primary 13Dxx, 16Dxx, 68W30;
 Secondary 16S37, 55Txx

{\bf Keywords.} local ring, Koszul complex, Yoneda Ext-algebra,

\noindent  Hilbert series, Macaulay2, Betti series, Avramov spectral sequence

\mysec{0. Introduction}
The aim of the present paper is to prove the following (if $V$ is a vector space over a field $k$ we  write $|V|$ instead of
${\rm dim}_kV$):

THEOREM 1 -- {\it Let $k$ be a field and $k[x,y,z,u]$ the polynomial ring in four variables of degree one and let
$$
R={k[x,y,z,u]\over (x^3,x^2y,(x+u)(y^2+z^2)+yzu,zu^2,u^3)} \leqno(1)
$$
be the graded quotient ring. Let
$$
P_R(Z)=\sum_{i \geq 0}|{\rm Tor}_i^R(k,k)|Z^i
$$
be the corresponding series of Betti numbers of $R$. Then $P_R(Z)$ is a transcendental function:
$$
 P_R(Z)= {Z(1+Z)^2 \over
(1+Z)(1-Z-Z^2)^2 \prod_{n=1}^\infty(1-Z^{2n})-1+2Z+2Z^2-4Z^3-4Z^4}
$$
}
If one replaces $k[x,y,z,u]$ with the formal power series ring $k[[x,y,z,u]]$ one obtains a local ring with
the same homological properties.
The Theorem might be of interest in itself (previously known examples of trancendental series
needed $\geq 5$ variables) but
more interesting are probably the methods used to obtain the example (described in section 1 below) and the
methods of proof described in section 2. There is now a possibility to obtain a new good insight in the
theory of the homology of local rings of embedding dimension 4. In section 3 we will in particular indicate
that the following ring which is a simpler variant of (1)
$$
R={k[x,y,z,u]\over (x^3,x^2y,(x+u)(y^2+z^2),zu^2,u^3)} \leqno(2)
$$
has a $P_R(Z)=(1-Z^2)(1+Z)/(1-3Z-3Z^3-2Z^5)$ i.e. a rational function {\it but}
 is such that the Yoneda Ext-algebra ${\rm Ext}^*_R(k,k)$ is
{\it not} finitely generated as an algebra. 
This last  phenomenon was also previously only known in the embedding dimension $\geq 5$ cases. We will also deduce some
other transcendental results and also show that there are quotients of $k[x,y,z,u]$ with an ideal with only
 {\it six} cubic generators which is non-Golod but has the multiplication in Koszul homology equal to zero
(the first examples of this last phenomenon with more relations - monomials of higher degrees and more variables -
were obtained by Lukas Katth\"an [KAT]).
 
\bigskip
\mysec{1. How the example in Theorem 1 was found}
In [A] (announced in 1980 [A-CR]), David Anick published the first example (example 7.1 page 29 of [A]) of an $R$ in 5 variables and 7 quadratic relations for which $P_R(Z)$ was transcendental.
 Here is a variant of his example with only 5 quadratic
relations which has the same property and which will be useful for us:
$$
R={k[x_1,x_2,x_3,x_4,x_5]\over (x_1^2,x_1x_2,x_1x_3+x_2x_4+x_3x_5,x_4x_5,x_5^2)} \leqno(3)
$$
In this case the Hilbert series of R is $H(t)=(1+3t+t^2-2t^3+t^4)/(1-t)^2$ and the Hilbert
series of the Koszul dual $R^!$ of R is
$$
R^!(t)={(1+t)^2\over(1-t-t^2)^2}\prod_{n=1}^\infty{1+t^{2n-1}\over 1-t^{2n}}
$$
and
$$
1/P_R(Z)=(1+1/Z)/R^!(Z)-H(-Z)/Z \leqno(4)
$$
Remark: For the case studied by Anick [A] we had that $H(t)=1+5t+7t^2$ and 
$$
R^!(t)={1\over (1-2t)^2}\prod_{n=1}^\infty{1+t^{2n-1}\over 1-t^{2n}}
$$ 
and the formula (4) was also valid in his case.

We now address the problem of constructing a 4-variable version of the ring (3).
We start with the ring $S_2=k[x,y,z,u]/(x^2,xy,zu,u^2)$ where $x,y,z,u$  correspond
to $x_1,x_2,x_4,x_5$ and try to add an extra relation corresponding to $x_1x_3+x_2x_4+x_3x_5$.
How to do this is not at all evident and instead we try to add a new relation which is
a linear combination with coefficients 0 or 1 of the 6 nonzero quadratic monomials
$yu,xu,z^2,yz,xz,y^2$ in $S_2$. There are 64 such linear combinations but each of them leads
to a rational $P_R(Z)$ (they all occur in the Appendix to [R4]).
Instead we turn to the study of cubic relations, i.e. we start with a new ring
$$
S={k[x,y,z,u]\over (x^3,x^2y,zu^2,u^3)}
$$
and try to add a linear combination (coefficients 0 or 1) of the 16 non-zero cubic monomials
$$
yu^2 , xu^2 , z^2u, yzu, xzu, y^2u, xyu, x^2u, z^3 , yz^2 , xz^2 , y^2z, xyz, x^2z, y^3 , xy^2
$$
in S. There are $2^{16}=65536$ cases to study and this can be done automatically using the following
input programme to Macaulay 2 which was written at my request several years ago by
Mike Stillman (it is an elegant version
of a programme I wrote for Macaulay 1 a long time ago (cf. [R3, pp. 294-296])

{\tts
binaries = (n) -> (

     if n === 0 then $\{\}$
     
     else if n === 1 then $\{\{0\},\{1\}\}$
     
     else (
     
     r := binaries(n-1);
     
     join(apply(r, i->prepend(0,i)),
    
        apply(r, i->prepend(1,i)))))

doit = (i) -> (

     h := hypers$\_{}$(0,i);

     J1 := J + ideal(h);

     << newline << flush;

     << ``--- n = `` << i << `` ideal = `` << hypers$\_{}$(0,i) << `` ---'' << flush;

     << newline << flush;

     E := res J1;
 
    << newline << flush;

     << ``  `` << betti E << flush;

     << newline << flush;

    A := (ring J1)/J1;
     
     C := res(coker vars A, LengthLimit=>6);
     
     << newline << flush;
     
     << ``  `` << betti C << flush;
      
     << newline << flush;

     << ``     ``    <<hilbertSeries(A,Order=>12)<< flush;

     << newline << flush; 

)

makeHyperplanes = (J) -> (

    I := matrix basis(3,coker gens J);

    c := numgens source I;

    m := transpose matrix binaries c;

    I * m)

R = QQ[x,y,z,u]

J = ideal(x\^{}3,x\^{}2*y,z*u\^{}2,u\^{}3)

time hypers = makeHyperplanes J;

time scan(numgens source hypers, i -> doit i);
}

\medskip
We obtain 20 possibilities for the $P_R(z)$ up to degree 6 presented here in increasing order (of course there can be
variations inside each case $P_i(z)$ due e.g. to different Koszul homology) :

$$
{\scriptstyle P_1(z)=1+4z+11z^2+27z^3+62z^4+137z^5+295z^6+\ldots \quad P_2(z)=1+4z+11z^2+27z^3+64z^4+152z^5+364z^6+\ldots}
$$
$$
{\scriptstyle P_3(z)=1+4z+11z^2+28z^3+69z^4+168z^5+407z^6+\ldots \quad P_4(z)=1+4z+11z^2+28z^3+69z^4+169z^5+414z^6+\ldots}
$$
$$
{\scriptstyle P_5(z)=1+4z+11z^2+28z^3+69z^4+172z^5+431z^6+\ldots \quad P_6(z)=1+4z+11z^2+28z^3+70z^4+173z^5+427z^6+\ldots}
$$
$$
{\scriptstyle P_7(z)=1+4z+11z^2+28z^3+70z^4+177z^5+451z^6+\ldots \quad P_8(z)=1+4z+11z^2+29z^3+74z^4+188z^5+476z^6+\ldots}
$$
$$
{\scriptstyle P_9(z)=1+4z+11z^2+29z^3+75z^4+193z^5+496z^6+\ldots \quad P_{10}(z)=1+4z+11z^2+29z^3+75z^4+193z^5+498z^6+\ldots}
$$
$$
{\scriptstyle P_{11}(z)=1+4z+11z^2+29z^3+75z^4+194z^5+503z^6+\ldots \quad P_{12}(z)=1+4z+11z^2+30z^3+80z^4+213z^5+567z^6+\ldots}
$$
$$
{\scriptstyle P_{13}(z)=1+4z+11z^2+30z^3+81z^4+218z^5+587z^6+\ldots \quad P_{14}(z)=1+4z+11z^2+30z^3+82z^4+224z^5+612z^6+\ldots}
$$
$$
{\scriptstyle P_{15}(z)=1+4z+11z^2+31z^3+86z^4+238z^5+660z^6+\ldots \quad P_{16}(z)=1+4z+11z^2+31z^3+87z^4+244z^5+685z^6+\ldots}
$$
$$
{\scriptstyle P_{17}(z)=1+4z+11z^2+31z^3+88z^4+249z^5+705z^6+\ldots \quad P_{18}(z)=1+4z+11z^2+32z^3+92z^4+263z^5+755z^6+\ldots}
$$
$$
{\scriptstyle P_{19}(z)=1+4z+11z^2+32z^3+93z^4+269z^5+780z^6+\ldots \quad P_{20}(z)=1+4z+11z^2+33z^3+99z^4+294z^5+877z^6+\ldots}
$$
Which one of these cases could give a transcendental $P_R(Z)$ ?

One can guess a possibility by studying how the transcendental $P_R(z)$ comes up in five variables
when we add linear combinations of quadratic monomials to the ring in 5 variables:
${\bf Q}[x,y,z,u,v]/(x^2,xy,uv,v^2)$ i.e. 
modifying the previous input programme by replacing R by $QQ[x,y,z,u,v]$, J by
$(x^2,x*y,u*v,v^2)$ and the line 

{\tt I := matrix basis(3,coker gens J);}

\noindent by the line

{\tt I := matrix basis(2,coker gens J);}

There are now ``only'' $2^{11}=2048$  cases to study and now we  
obtain only 10 possibilities for the $P_R(z)$ up to degree 6 presented here in increasing order:
$$
{\scriptstyle P_{1}(z)=1+5z+15z^2+37z^3+82z^4+170z^5+337z^6+\ldots \quad P_{2}(z)=1+5z+15z^2+38z^3+88*z^4+192z^5+406z^6+\ldots}
$$
$$
{\scriptstyle P_{3}(z)=1+5z+15z^2+38z^3+89z^4+199z^5+432z^6+\ldots \quad P_{4}(z)=1+5z+15z^2+38z^3+91z^4+216z^5+516z^6+\ldots}
$$
$$
{\scriptstyle P_{5}(z)=1+5z+15z^2+39z^3+95z^4+222z^5+506z^6+\ldots \quad P_{6}(z)=1+5z+15z^2+39z^3+96z^4+231z^5+553z^6+\ldots}
$$
$$
{\scriptstyle P_{7}(z)=1+5z+15z^2+39z^3+97z^4+237z^5+575z^6+\ldots \quad P_{8}(z)=1+5z+15z^2+39z^3+99z^4+254z^5+659z^6+\ldots}
$$
$$
{\scriptstyle P_{9}(z)=1+5z+15z^2+40z^3+104z^4+268z^5+689z^6+\ldots \quad P_{20}(z)=1+5z+15z^2+41z^3+112z^4+306z^5+836z^6+\ldots}
$$
There are e.g. 1024 cases of an extra quadratic relation corresponding to $P_1(z)$ and 576 cases corresponding to $P_4(z)$, but
only 8 for $P_2(z)$ including the case of the extra relation $xz+yu+zv$ leading to the transcendental series mentioned
in the introduction (the other 7 cases of $P_2(z)$ lead to the same total series).
Furthermore, there are also 8 cases for $P_6(z)$ and they all lead to cases where the Yoneda Ext-algebra 
${\rm Ext}^*_R(k,k)$ is not finitely generated but the series of Betti numbers is rational 
(this phenomenon was first found by me in [R1] (inspired by Lemaire [LEM]).

We now try to use this information to try to guess what happens for the 65536 embedding dimension 4 cases above.
In this situation there are only 80 cases of $P_{10}(z)$ and only two of them correspond to adding a cubic relation
with is a sum of 5 cubic monomials (the other cases need more monomials):
case 9257 which corresponds to $xy^2+xz^2+y^2u+yzu+z^2u$ and case 13345 which corresponds to
$xy^2+xyz+xz^2+y^2u+z^2u$. Here one of these cases can be  obtained from the other one by permuting x and u.
We can therefore restrict ourselves to the homological study of the ring 
$$
k[x,y,z,u]\over(x^3,x^2y,(x+u)(y^2+z^2)+yzu,zu^2,u^3) \leqno(1)
$$  
mentioned in the introduction and which will be the subject of the next section.    
\bigskip
\mysec{2. The treatment of the example in Theorem 1}
We will now use the programme {\tt DGAlgebra} of [MOO] on the ring (1). 

The Koszul homology of the ring has the Betti numbers given by Macaulay:
$$
\vbox{\halign{
 #\quad&#\quad&#\quad&#\quad&#\quad&#\quad\cr
total:& 1& 5& 12& 12& 4\cr
0:& 1& .& .&  .&   .\cr
1:& .& .& .&  .&  .\cr
2:& .& 5& 2&  .&  .\cr
3:& .& .& 1&  .& .\cr
4:& .& .& 8& 10& 3\cr
5:& .& .& .&  .& .\cr
6:& .& .& 1&  2& 1\cr
}}
$$
From this we can guess that the multiplication in $HKR={\rm Tor}_*^{k[x,y,z,u]}(k[x,y,z,u]/J,k)$,
where $J$ is the ideal in (1)
is rather complicated. We determine details about that multiplication as a
preparation for using the Avramov spectral sequence ([AV0],[AV1]) where $\tilde R = k[x,y,z,u]$ and $R={\tilde R}/J$:
$$
E^2_{p,q}={\rm Tor}^{HKR}_{p,q}(k,k) => {\rm Tor}^R_*(k,k)//{\rm Tor}^{\tilde R}_*(k,k) \leqno(5)
$$
by using the following infile for Macaulay:

{\tt
loadPackage(``DGAlgebras'')

R=QQ[x,y,z,u]/ideal(x\^{}3,x\^{}2*y,z*u\^{}2,u\^{}3,x*y\^{}2+x*z\^{}2+y\^{}2*u+y*z*u+z\^{}2*u)

res(ideal R); betti oo

res(coker vars R,LengthLimit => 6); betti oo

HKR=HH koszulComplexDGA(R)

generators HKR

for n from 1 to length(generators HKR) list degree X{\_}n

ideal HKR

I=ideal(vars HKR)

I\^{}2; trim(oo)

I\^{}3; trim(oo)

res(coker vars HKR,LengthLimit => 3);

betti oo

res(coker vars HKR,LengthLimit => 4);

betti oo
}

From the Macaulay output of this input file we see that the augmentation ideal
$I$ of $HKR$ satisfies $I^3=0$, and that $I^2$ is generated by
$$
{\scriptstyle X_4X_5,X_2X_5,X_1X_5,X_3X_4,X_2X_4,X_1X_4,X_2X_3,X_4X_7,X_2X_7,X_1X_7,X_5X_6,X_4X_6,X_3X_6.X_6X_7}
$$
 Therefore we will be able to calculate the left side
of the Avramov spectral sequence (5) exactly as in [R4], where we used the Theorem 1.3 on page 310
of Clas L\"ofwall's paper [L1], provided we can determine the Hilbert series
$HKR^!(z,z)$
of the Koszul dual $HKR^!$ of $HKR$ and the Hilbert series $HKR(-z,z)$ of $HKR$.
From the output file we see that $HKR$ can be presented as follows:
It has five generators $X_1,X_2,X_3,X_4,X_5$ of degree (1,3),two generators $X_6,X_7$ of degree
(2,4), three generators $X_8,X_9,X_{10}$ of degrees (2,5),(2,6),(2,8) respectively, four generators
$X_{11},X_{12},X_{13},X_{14}$ of degree (3,7), two generators $X_{15},X_{16}$ of degree (3,9),
two generators $X_{17},X_{18}$ of degree (4,8) and one generator $X_{19}$ of degree (4,10). 

Furthermore the first five variables $X_1,X_2,X_3,X_4,X_5$ skew-commute, the variables $X_6,X_7,X_8,X_9,X_{10}$
commute among themselves and with $X_1,X_2,X_3,X_4,X_5$, etc. 

The relations are as follows: The product of the last 12 variables $X_8,\ldots X_{19}$
with all the 19 variables $X_1 \ldots X_{19}$ are 0.
For the first five generators $X_1,X_2,X_3,X_4,X_5$ we only have the following three quadratic relations:
$X_3X_5,X_1X_3+X_2X_4-X_4X_5,X_1X_2$, and if we also take into account the additional variables $X_6,X_7$
 we have the extra quadratic relations $X_1X_6,X_2X_6,X_3X_7,X_5X_7,X_6X_6,X_7X_7$.
It follows that the if 
$$
C={k[X_1,X_2,X_3,X_4,X_5,X_6,X_7]\over (X_3X_5,X_1X_3+X_2X_4-X_4X_5,X_1X_2,X_1X_6,X_2X_6,X_3X_7,X_5X_7,X_6X_6,X_7X_7)} \leqno(6)
$$
then the Koszul dual  $HKR^!$ of $HKR$ can be presented as a coproduct 
$$
C^! \sqcup k<x_8,x_9,x_{10},x_{11},x_{12},x_{13},x_{14},x_{15},x_{16},x_{17},x_{18},x_{19}> \leqno(7)
$$
where the second algebra in the coproduct is the free algebra on the dual generators $x_8,\ldots x_{19}$ of $X_8,\ldots X_{19}$.

Furthermore, it is well-known (cf. e.g. [LEM]) that for the coproduct $A \sqcup B$ of two graded connected algebras,
we have the following formula for the respective Hilbert series
$$
{1\over {(A \sqcup B)}(z)} +1 = {1\over A(z)}+{1\over B(z)} 
$$
(this works also for  the bigraded case) and it follows that 
 $$
{1\over HKR^!(z,z)}+1={1\over C^!(z,z)}+1-3z^3-6z^4-3z^5 \leqno(8)
$$ 
and we have therefore reduced our problem to calculate $C^!(z,z)$ where C only involves the first 7 variables
 $X_1,X_2,X_3,X_4,X_5,X_6,X_7$
 
We start with an algebra $D$ that involves only the first 5 variables $X_1,X_2,X_3,X_4,X_5$:
$$
D={k[X_1,X_2,X_3,X_4,X_5]\over (X_3X_5,X_1X_3+X_2X_4-X_4X_5,X_1X_2)}
$$
The variables skewcommute. The Koszul dual $D^!$ is the quotient of the free algebra in the dual variables $x_i$ of the
$X_i$ :
$$
D^!={k<x_1,x_2,x_3,x_4,x_5>\over([x_1,x_4],[x_1,x_5],[x_2,x_3],[x_2,x_5],[x_3,x_4],
[x_1,x_3]+[x_4,x_5],[x_2,x_4]+[x_4,x_5]}
$$ 
where the $[x_i,x_j]$ are Lie commutators.
From this one sees exactly as in L\"ofwall's and mine study of the Anick example [LR1] that the underlying Lie algebra $g_D$ of $D^!$
sits in the middle of an extension
$$
0 \longrightarrow \oplus_{i=1}^{\infty}a_i  \longrightarrow g_D  \longrightarrow f_1 \times f_2 \longrightarrow 0
$$
i.e. is the extension of the product of two free Lie algebras $f_1$ generated by $x_1$ and $x_2$ and
$f_2$ generated by $x_3$ and $x_5$  with the infinite abelian Lie algebra
$$
a=\oplus_{i=1}^{\infty}a_i
$$
where the $a_i$ are onedimensional generated by $a_1=x_4$, $a_2=[x_5,x_4]$, $a_3=[x_5,[x_5,x_4]]$
 etc. which commute and  by the 2-cocycle $\gamma: f_1 \times f_2 \rightarrow a$
defined by $\gamma(x_1,x_3) = -\gamma(x_3,x_1)=a_2$ and  $\gamma(x_2,x_5) = -\gamma(x_5,x_2)=a_2$.
It follows that 
$$
D^!(z)={1\over (1-2z)^2}\prod_{n=1}^\infty{1\over(1-z^n)}  \leqno(7)
$$
But since we need $D^!(z,z)$ we have to replace the $z$ by $z^2$ in the formula (9) to get $D^!(z,z)$.
We next want to incorporate the variables $X_6,X_7$ (which {\it commute} with the $X_1,X_2,X_3,X_4,X_5$)
into the picture. It turns out that we have to study the underlying Lie algebra of the Koszul dual $C^!$ of
the $C$ in (6), i.e. the quotient of 
$$
k<x_1,x_2,x_3,x_4,x_5,x_6,x_7>
$$
 with the ideal generated by
$$
[x_1,x_4],[x_1,x_5],[x_2,x_3],[x_2,x_5],[x_3,x_4],[x_1,x_3]+[x_4,x_5],[x_2,x_4]+[x_4,x_5],
$$
as before and the extra ideal generators:
$$
[x_1,x_7],[x_2,x_7],[x_4,x_7],[x_6,x_7],[x_3,x_6],[x_4,x_6],[x_6,x_7]
$$
where $x_6$ and $x_7$ are dual to $X_6$ and $X_7$.
From the preceding presentation one sees that the Lie algebra $g_C$ of $C^!$ sits in the middle of
an extension:
$$
0 \longrightarrow \oplus_{i=1}^{\infty}a_i  \longrightarrow g_C  \longrightarrow F_1 \times F_2 \longrightarrow 0 \leqno(10)
$$
where $F_1$ is the free Lie algebra generated by $x_1,x_2,x_6$ of degrees $1,1,2$ and 
$F_2$ is the free Lie algebra generated by $x_3,x_5,x_7$ of degrees $1,1,2$ and the $a_i$:s are as before,
and the $x_6$ and $x_7$ operate in the trivial way on the $a_i$.
Furthermore we have a new cocycle
$$
\Gamma : F_1 \times F_2 \longrightarrow  \oplus_{i=1}^{\infty}a_i
$$ 
where as before $\Gamma$ is defined by its value on the generators by $-\Gamma(x_3,x_1)=\Gamma(x_1,x_3)=a_2$ 
and $-\Gamma(x_5,x_2)=\Gamma(x_2,x_5)=a_2$ and being zero for all other pairs of generators.
By calculating $H^2(F_1 \times F_2,a)$ for (8) one sees that the relations are those of $C^!$ and therefore
the Hilbert series
$$
C^!(z)={1\over (1-2z-z^2)^2}\prod_{n=1}^\infty{1\over(1-z^n)}  \leqno(11) 
$$
if the variables $x_1,x_2,x_3,x_4,x_5,x_6,x_7$ are given the degrees $1,1,1,1,1,2,2$. But for the calculation of
the series of $HKR^!(z,z)$ these variables should be given the degrees 

\noindent $2,2,2,2,2,3,3$ so that we get the series
$$
{1\over (1-2z^2-z^3)^2}\prod_{n=1}^\infty{1\over(1-z^{2n})}  \leqno(12)
$$
Using the formulae (8) and (12) we obtain 
$$
{1 \over HKR^!(z,z)}=(1-2z^2-z^3)^2\prod_{n=1}^\infty(1-z^{2n})-3z^3-6z^4-3z^5 \leqno(13) 
$$
But $I^3=0$ and
$$
{\scriptstyle I^2=(X_4X_5,X_2X_5,X_1X_5,X_3X_4,X_2X_4,X_1X_4,X_2X_3,X_4X_7,X_2X_7,X_1X_7,X_5X_6,X_4X_6,X_3X_6.X_6X_7)}
$$ 
and this gives that $HKR(x,y)=1+5xy+5xy^2+6xy^3+3xy^4+7x^2y^2+6x^2y^3+x^2y^4$, so that
$$
 HKR(-z,z)=1-5z^2-5z^3+z^4+3z^5+z^6 \leqno(14)
$$
and finally the formula
$$
{1\over P_{HKR}(z,z)}=(1+1/z)/HKR^!(z,z)-HKR(-z,z)/z
$$
gives that 
$$
{1\over P_{HKR}(z,z)}=(1+1/z)[(1-2z^2-z^3)^2\prod_{n=1}^{\infty}(1-z^{2n})-3z^3-6z^4-3z^5)]
$$
$$
-(1-5z^2-5z^3+z^4++3z^5+z^6)/z
$$
and {\it if the Avramov spectral sequence degenerated} this series should be the same as 

\noindent $(1+z)^4/P_R(z)$ so the final result should be
$$
 P_R(z)= {z(1+z)^2 \over
(1+z)(1-z-z^2)^2 \prod_{n=1}^\infty(1-z^{2n})-1+2z+2z^2-4z^3-4z^4}
$$
$$
=1+4z+11z^2+29z^3+75z^4+193z^5+498z^6+1289z^7+3341z^8+8663z^9+22466z^{10}+\ldots
$$
as asserted in our Theorem 1. The following calculation in Macaulay2

{\tt res(coker vars R, LengthLimit => 11); betti oo}

\noindent gives the diagram of Betti numbers (16) below whose first row (the total
Betti numbers) gives support to this assertion.
But the Theorem 5.9 of [AV0] gives first that the differentials $d^r$ of the
Avramov spectral sequence are 0 for $r>2$ but furthermore since the dimensions of the two sides of the Avramov spectral
sequence are the same in degrees $\leq 10$ it follows from the proof of Theorem 5.9 in [AV0] and the structure of
the matrix Massey products that there are no non-zero differentials, and Theorem 1 is proved.
We can also get a bigraded more precise version of all this:
Indeed, we also have a series 
$$
P_R(x,y)=\sum_{p,q}|{\rm Tor}_{p,q}^R(k,k)|x^py^q
$$
which takes into account the extra grading of $R$ (we have that $P_R(z,1)$ is the old $P_R(z)$)
given as 
$$
P_R(x,y)={(1+xy)^4\over(1+1/x)/A-H/x} \leqno(15)
$$
where
$$
  1/A=(1-2x^2y^3-x^3y^4)^2\prod_{n=1}^{\infty}(1-x^{2n}y^{3n})-x^3y^5-x^3y^6-x^3y^8-4x^4y^7-2x^4y^9-2x^5y^8-x^5y^{10}
$$
and $H$ is equal to
$$
1-5x^2y^3-2x^3y^4-x^3y^5-x^3y^6-x^3y^8-4x^4y^7-2x^4y^9-2x^5y^8-x^5y^{10}+7x^4y^6+6x^5y^7+x^6y^8
$$
and the expansion of $P_R(x,y)$ of (15) up to degree 11 gives the same diagram of graded Betti numbers that was 
earlier found by the Macaulay2 calculation:

$$
\vbox{\halign{
 #\quad\hfill&#\quad\hfill&#\quad\hfill&#\quad\hfill&#\quad\hfill&#\quad\hfill&#\quad\hfill&#\quad\hfill&#\quad\hfill&#\quad\hfill&#\quad\hfill&#\quad\hfill&#\quad\hfill\cr
total:&1&4&11&29&75&193&498&1289&3341&8663&22466&58264\cr
    0:&1&4& 6& 4& 1&  .&  .&   .&   .&   .&   .&.\cr
1:&.&.& 5& 22& 38&  32&  13&   2&   .&   .&   .&.\cr
2:&.&.& .& 1& 22&  92& 171&  169&   92&  26& 3&.\cr
3:&.&.& .& 1& 8&  34&  135&  389&   689& 744& 493&196\cr
4:&.&.& .& .& .&  10&  87&   346&   923& 1870& 2835&3066\cr
5:&.&.& .& 1& 6&  15&  22&   96&    645& 2426& 5739&9705\cr
6:&.&.& .& .& .&  10&  65&  186&    348&  964& 4302&14387\cr
7:&.&.& .& .& .&  .&   2&   73&     459&  1446& 3296&8595\cr
8:&.&.& .& .& .&  .&   2&   20&     112&   669& 3025&9538\cr
9:&.&.& .& .& .&  .&  .&     .&      30&   320& 1668&6469\cr
10:&.&.&.& .& .&  .&   1&    8&      28&   62&   396&3282\cr
11:&.&.& .& .& .&  .&  .&   .&       15&   129&  510&1377\cr
12:&.&.& .& .& .&  .&  .&   .&        .&   3&    153&1278\cr
13:&.&.& .& .& .&  .&  .&   .&        .&   3&     36&246\cr
14:&.&.& .& .& .&  .&  .&   .&    .&   .&   .&60\cr
15:&.&.& .& .& .&  .&  .&   .&   .&   1&   10&45\cr
16:&.&.& .& .& .&  .&  .&   .&   .&   .&   .&20\cr 
}} \leqno(16)
$$

\mysec{3. Other embedding dimension 4 cases}
\medskip
The example we have just studied is certainly not unique. Let us illustrate this with two examples:
If one adds the relation $z^3$ to the example in Theorem 1
we obtain the following
$$
R_a={k[x,y,z,u]\over (x^3,x^2y,(x+u)(y^2+z^2)+yzu,z^3,zu^2,u^3)}
$$
which has Hilbert series
$$
1+3T+6T^2+4T^3-T^4-7T^5 \over 1-T
$$
Furthermore the homology of the Koszul complex of $R_a$ is
$$
\vbox{\halign{
 #\quad&#\quad&#\quad&#\quad&#\quad&#\quad\cr
total:& 1& 6& 18& 20& 7\cr
0:& 1& .& .&  .&   .\cr
1:& .& .& .&  .&  .\cr
2:& .& 6& 2&  .&  .\cr
3:& .& .& 2&  .& .\cr
4:& .& .&14& 20& 7\cr
}}
$$
We get same Hilbert series and the same diagram of Koszul homology if we replace $z^3$ with $y^3+z^3$ to get the ring
$$
R_b={k[x,y,z,u]\over (x^3,x^2y,(x+u)(y^2+z^2)+yzu,y^3+z^3,zu^2,u^3)}
$$
Both these ring have transcendental series of Betti numbers. These two series can be determined and they are different. For the case of $R_a$ we have
using {\tt DGAlgebras} that $HKR_a$ is generated by 6 skew-commuting 
variables $X_1,X_2,X_3,X_4,X_5,X_6$ of degree (1,3) 
(and 27 variables of higher degrees). The first 6  variables satisfy the
following relations $X_5X_6,X_3X_6,X_3X_5,X_1X_3+X_2X_4-X_4X_6,X_1X_2$ and
the Hilbert series of the corresponding Koszul dual is 
$$
1\over (1-2z)(1-3z)\prod_{n=1}^{\infty}(1-z^n)
$$
On the other hand, using {\tt DGAlgebras} for $R_b$ we find that $HKR_b$
is still generated by 6 variables $X_1,X_2,X_3,X_4,X_5,X_6$ of degree (1,3)
(but only 24 variables of higher degrees). In this case the first 6  variables satisfy the following relations $X_3X_6,X_1X_3+X_2X_4-X_4X_6,X_1X_2$
and in  this case the Hilbert series of the corresponding Koszul dual is
$$
1\over (1-2z)^2 (1-z)\prod_{n=1}^{\infty}(1-z^n)
$$
Now we can continue our analysis as in section 2 and we
obtain the following formulae:
$$
P_{R_a}(Z)= {Z(1+Z)^2 \over
(1-Z-Z^2)(1-3Z^2-Z^3) \prod_{n=1}^\infty(1-Z^{2n})-1+2Z+3Z^2-6Z^3-7Z^4}
$$
and
$$
 P_{R_b}(Z)= {Z(1+Z)^2 \over
(1-Z)(1-2Z^2-Z^3)^2 \prod_{n=1}^\infty(1-Z^{2n})-1+2Z+3Z^2-6Z^3-7Z^4}
$$
Note also that we could have started with many more examples in embedding dimension 4 than just our
 variation of the Anick case that we took in section 1.

Now let us also briefly analyze the case when
$$
S={k[x,y,z,u]\over (x^3,x^2y,(x+u)(y^2+z^2),zu^2,u^3)}
$$
Using {\tt DGAlgebras} we obtain as before that  
$$
P_S(Z)={(1-Z)(1+z)^2 \over 1-3Z+3Z^3-2Z^5}
$$
but we can also prove that neither ${\rm Ext}^*_S(k,k)$ nor ${\rm Ext}^*_{HKS}(k,k)$  
are finitely generated as algebras equipped with the Yoneda product.
The case $S$ is the 4-variable of the 5-variable $k[x,y,z,u,v]/(x^2,xy,xz+zv,uv,v^2)$,
a variant of which we treated in [R1], inspired by Lemaire [LEM].
An alternative way to treat $S$ is to observe that the quotient map
$$
{k[x,y,z,u]\over (x^3,x^2y,zu^2,u^3)} \longrightarrow S
$$
is a Golod map [LEV] so that $k[x,y,z,u]\rightarrow S$ is a composite of three Golod maps,
and use [R2].
 
All this leads support to the surmise that the embedding dimension 4 case could be equally
complicated as the general case, and one could even pose the

QUESTION 1: Let ${\cal E}4$ be the set of series $\sum_{n\geq 0}|{\rm Tor}^R_n(k,k)|Z^n$ for 
$(R,m)$ a local commutative noetherian ring of embedding dimension $\leq 4$ where $k=R/m$. Can ${\cal E}4$ 
be added to the list of 17 series of [A-Gu] that are all rationally related ?

One would certainly get a smaller set than ${\cal E}4$ if one restricted oneself to rings of the form
$k[x,y,z,u]/(f_1,\ldots f_s)$ where the relations $f_i$ were of degree $\leq 3$ or $\leq 4$ .

So far no embedding dimension 4 variant of the rings in [LR2] has been found.
But even in this restricted case of degree 3 relations one can get surprizes:
The ring 
$$
S_I=Q[x,y,z,u]/(u^3,xy^2,(x+y)z^2,x^2u+zu^2,y^*u+xzu,y^2z+yz^2)
$$
constructed in [R5] (inspired by Lukas Katth\"an [KAT])
which has only six relations of degree 3, has as homology of the Koszul complex:
$$
\vbox{\halign{
 #\quad&#\quad&#\quad&#\quad&#\quad&#\quad\cr
Total: & 1& 6& 12& 9& 2\cr
0: & 1& .& .&  .&   .\cr
1: & .& .& .&  .&  .\cr
2: & .& 6& .&  .&  .\cr
3: & .& .& 12&  5& 1\cr
4: & .& .& .&  4& .\cr
5: & .& .&  .& .& 1\cr
}}
$$
Furthermore the multiplication of elements of positive degree in this homology is zero (use {\tt DGAlgebra}), and
$$
{(1+z)^4\over(1-6z^2-12z^3-9z^4-2z^5)}=1+4z+12z^2+40z^3+130z^4+422z^5+1376z^6+4476z^7+\ldots
$$
whereas the series of the total Betti numbers of $S_I$ is given by Macaulay2 as:
$$
1+ 4z+ 12z^2+ 40z^3+ 130z^4+ 421z^5+ 1371z^6+ 4454z^7+\ldots
$$
so that $S_I$ is {\it not} a Golod ring.
In [R5] we show that if one accepts relations of degree 3 {\it and} 4, there are 
hundreds of such exotic non-Golod rings.

QUESTION 2: For local rings $R$ of embedding dimension 4 with only cubic relations,
are my $R$ in Theorem 1 and the $R_a$ and $R_b$ essentially the only examples where we have
a transcendental $P_R(Z)$ ?
 
\mysec{4. References}
\medskip
{\frenchspacing\raggedbottom
\def\myref#1{\smallskip\item{[{\bf #1}]}}

\myref{A} D.J.  Anick, {\it A counterexample to a conjecture of Serre},
Ann. of Math. {\bf 115}, 1982, pp. 1-33; {\it Correction}, Ann. of Math. {\bf 116} ,
p. 661.

\myref{A-CR} D.J. Anick, {\it Construction d'espaces de lacets et d'anneaux locaux \`a s\'eries de Poincaré-Betti non rationnelles},
C.R. Acad. Sc. Paris {\bf 290}, (1980), pp. A733-A736.

\myref{A-Gu} D.J.  Anick and Tor H. Gulliksen, {\it Rational dependence among Hilbert and Poincare series},
Journal of Pure and Applied Algebra {\bf 38}, 1985, pp. 135-157.

\myref {AV0} Avramov, L. L.,{\it The Hopf algebra of a local ring. (Russian)},
Izv. Akad. Nauk SSSR Ser. Mat. {\bf 38}, (1974), pp. 253-277. English translation [Math. USSR-Izv.{\bf  8} (1974), 259-284].

\myref {AV1} Avramov, L. L.,{\it Obstructions to the existence of multiplicative structures on minimal free resolutions.}
 Amer. J. Math., {\bf 103} (1981), pp. 1-31.

\myref{KAT} L. Katth\"an, {\it A non-Golod ring with a trivial product on its Koszul homology},

{\tt http://arxiv.org/pdf/1511.04883.pdf}

\myref{LEM} Lemaire, Jean-Michel,{\it Alg\`ebres connexes et homologie des espaces de lacets}, Lecture Notes in Mathematics, {\bf  422}, Springer-Verlag, Berlin-New York, 1974.

\myref{LEV} Levin, Gerson, {\it Local rings and Golod homomorphisms}, Journal of Algebra {\bf 37}, 1975, pp 266–289.

\myref{L1} L\"ofwall, C., {\it On the subalgebra generated by one-dimensional elements in the Yoneda {\rm Ext}--algebra},
 in Algebra, algebraic topology, and their interactions,
(       J.--E. Roos, ed),       Lecture Notes in Math., vol.{\bf  1183}, Springer-Verlag,
Berlin--New York, 1986, pp. 291-338.

\myref{LR1} C. L\"ofwall and J.-E. Roos, {\it Cohomologie des alg\`ebres de Lie gradu\'ees et s\'eries de Poincar\'e-Betti
 non rationnelles}, C.R. Acad. Sc. Paris {\bf 290}, 1980, pp. A733-A736.

\myref{LR2} C. L\"ofwall and J.-E. Roos, {\it  A nonnilpotent 1-2-presented graded Hopf algebra whose Hilbert series converges in the \
unit circle},
 Adv. Math.{\bf  130}, 1997, pp. 161-200.

\myref{MOO} Moore, Frank, {\tt DGAlgebras.} A package for {\tt Macaulay2} 

{\tt http://www.math.uiuc.edu/Macaulay2/Packages/}

\myref{R1} J.-E. Roos, {\it Relations between the Poincar\'e-Betti series of Loop Spaces and of Local rings},  Springer Lecture Notes \
in Math.{\bf 740},1979, 285-322.
\myref{R2} J.-E. Roos, {\it On the use of graded Lie algebras in the theory of local rings}, Commutative algebra: Durham 1981 (R. Y. S\
harp, ed.) London Math.
Soc. Lecture Notes Ser. vol. 72, Cambridge Univ. Press, Cambridge, 1982, pp. 204–230.

\myref{R3} J.-E. Roos, {\it A computer-aided study of the graded Lie-algebra of a local commutative noetherian ring} (with an Appendix
by Clas L\"ofwall),  Journal of Pure and Applied Algebra {\bf 91 }, 1994, pp. 255-315.

\myref{R4} J.-E. Roos, {\it Homological properties of the homology algebra of the Koszul complex of a local ring: Examples and questions},  Journal of Algebra {\bf 465}, 2016, pp. 399-436.

\myref{R5} J.-E. Roos, {\it On some unexpected rings that are close to Golod rings}, in preparation.

} 
\par}
\end